\documentclass[12pt]{article}
\usepackage{amsmath,amssymb,amsthm,latexsym,amscd}

\title{Definable families of theories,
\\ related calculi and ranks\footnote{
\newline\indent \ \ \ This research was partially
supported by Committee of Science in Education and Science
Ministry of the Republic of Kazakhstan (Grants No. AP05132349,
AP05132546) and Russian Foundation for Basic Researches (Project
No. 17-01-00531-a). } }
\author{Nurlan D.
Markhabatov, Sergey V.
Sudoplatov\footnote{nur\underline{\,\,\,}24.08.93@mail.ru,
sudoplat@math.nsc.ru}}
\date{}
\begin{document}
\maketitle

\begin{abstract}
We consider sentence-definable and diagram-definable subfamilies
of given families of theories, calculi for these subfamilies, as
well dynamics and characteristics of these subfamilies with
respect to rank and degree.

\end{abstract}

The rank for families of theories was introduced and studied in
general context in \cite{RSrank}. All possible values of the ranks
and degrees for families of all theories in given languages were
described in \cite{MarkhSud}. In the present paper we consider
sentence-definable and diagram-definable subfamilies of given
families of theories, calculi for these subfamilies, as well as
dynamics and characteristics of these subfamilies with respect to
rank and degree.

The paper is organized as follows. In Section 1, preliminary
notions, notations and results are collected. In Section 2, we
consider calculi subfamilies of families of theories as well as
links for sentence-definable and diagram-definable subfamilies.
Compactness and $E$-closeness for definable subfamilies are
studied in Section 3. In Section 4 we consider dynamics of ranks
with respect to definable subfamilies of theories and prove the
existence of subfamilies of given rank.

\section{Preliminaries}

Throughout we consider families $\mathcal{T}$ of complete
first-order theories of a language $\Sigma=\Sigma(\mathcal{T})$.

Throughout the paper we consider complete first-order theories $T$
in predicate languages $\Sigma(T)$ and use the following
terminology in \cite{RSrank, cs, cl, lut, ccct, rest, lft, at}.

Let $P=(P_i)_{i\in I}$, be a family of nonempty unary predicates,
$(\mathcal{A}_i)_{i\in I}$ be a family of structures such that
$P_i$ is the universe of $\mathcal{A}_i$, $i\in I$, and the
symbols $P_i$ are disjoint with languages for the structures
$\mathcal{A}_j$, $j\in I$. The structure
$\mathcal{A}_P\rightleftharpoons\bigcup\limits_{i\in
I}\mathcal{A}_i$\index{$\mathcal{A}_P$} expanded by the predicates
$P_i$ is the {\em $P$-union}\index{$P$-union} of the structures
$\mathcal{A}_i$, and the operator mapping $(\mathcal{A}_i)_{i\in
I}$ to $\mathcal{A}_P$ is the {\em
$P$-operator}\index{$P$-operator}. The structure $\mathcal{A}_P$
is called the {\em $P$-combination}\index{$P$-combination} of the
structures $\mathcal{A}_i$ and denoted by ${\rm
Comb}_P(\mathcal{A}_i)_{i\in I}$\index{${\rm
Comb}_P(\mathcal{A}_i)_{i\in I}$} if
$\mathcal{A}_i=(\mathcal{A}_P\upharpoonright
A_i)\upharpoonright\Sigma(\mathcal{A}_i)$, $i\in I$. Structures
$\mathcal{A}'$, which are elementary equivalent to ${\rm
Comb}_P(\mathcal{A}_i)_{i\in I}$, will be also considered as
$P$-combinations.

Clearly, all structures $\mathcal{A}'\equiv {\rm
Comb}_P(\mathcal{A}_i)_{i\in I}$ are represented as unions of
their restrictions $\mathcal{A}'_i=(\mathcal{A}'\upharpoonright
P_i)\upharpoonright\Sigma(\mathcal{A}_i)$ if and only if the set
$p_\infty(x)=\{\neg P_i(x)\mid i\in I\}$ is inconsistent. If
$\mathcal{A}'\ne{\rm Comb}_P(\mathcal{A}'_i)_{i\in I}$, we write
$\mathcal{A}'={\rm Comb}_P(\mathcal{A}'_i)_{i\in
I\cup\{\infty\}}$, where
$\mathcal{A}'_\infty=\mathcal{A}'\upharpoonright
\bigcap\limits_{i\in I}\overline{P_i}$, maybe applying
Morleyzation. Moreover, we write ${\rm
Comb}_P(\mathcal{A}_i)_{i\in I\cup\{\infty\}}$\index{${\rm
Comb}_P(\mathcal{A}_i)_{i\in I\cup\{\infty\}}$} for ${\rm
Comb}_P(\mathcal{A}_i)_{i\in I}$ with the empty structure
$\mathcal{A}_\infty$.

Note that if all predicates $P_i$ are disjoint, a structure
$\mathcal{A}_P$ is a $P$-combination and a disjoint union of
structures $\mathcal{A}_i$. In this case the $P$-combination
$\mathcal{A}_P$ is called {\em
disjoint}.\index{$P$-combination!disjoint} Clearly, for any
disjoint $P$-combination $\mathcal{A}_P$, ${\rm
Th}(\mathcal{A}_P)={\rm Th}(\mathcal{A}'_P)$, where
$\mathcal{A}'_P$ is obtained from $\mathcal{A}_P$ replacing
$\mathcal{A}_i$ by pairwise disjoint
$\mathcal{A}'_i\equiv\mathcal{A}_i$, $i\in I$. Thus, in this case,
similar to structures the $P$-operator works for the theories
$T_i={\rm Th}(\mathcal{A}_i)$ producing the theory $T_P={\rm
Th}(\mathcal{A}_P)$\index{$T_P$}, being {\em
$P$-combination}\index{$P$-combination} of $T_i$, which is denoted
by ${\rm Comb}_P(T_i)_{i\in I}$.\index{${\rm Comb}_P(T_i)_{i\in
I}$}

Notice that $P$-combinations are represented by generalized
products of structures \cite{FV}.

For an equivalence relation $E$ replacing disjoint predicates
$P_i$ by $E$-classes we get the structure
$\mathcal{A}_E$\index{$\mathcal{A}_E$} being the {\em
$E$-union}\index{$E$-union} of the structures $\mathcal{A}_i$. In
this case the operator mapping $(\mathcal{A}_i)_{i\in I}$ to
$\mathcal{A}_E$ is the {\em $E$-operator}\index{$E$-operator}. The
structure $\mathcal{A}_E$ is also called the {\em
$E$-combination}\index{$E$-combination} of the structures
$\mathcal{A}_i$ and denoted by ${\rm Comb}_E(\mathcal{A}_i)_{i\in
I}$\index{${\rm Comb}_E(\mathcal{A}_i)_{i\in I}$}; here
$\mathcal{A}_i=(\mathcal{A}_E\upharpoonright
A_i)\upharpoonright\Sigma(\mathcal{A}_i)$, $i\in I$. Similar
above, structures $\mathcal{A}'$, which are elementary equivalent
to $\mathcal{A}_E$, are denoted by ${\rm
Comb}_E(\mathcal{A}'_j)_{j\in J}$, where $\mathcal{A}'_j$ are
restrictions of $\mathcal{A}'$ to its $E$-classes. The
$E$-operator works for the theories $T_i={\rm Th}(\mathcal{A}_i)$
producing the theory $T_E={\rm Th}(\mathcal{A}_E)$\index{$T_E$},
being {\em $E$-combination}\index{$E$-combination} of $T_i$, which
is denoted by ${\rm Comb}_E(T_i)_{i\in I}$\index{${\rm
Comb}_E(T_i)_{i\in I}$} or by ${\rm
Comb}_E(\mathcal{T})$\index{${\rm Comb}_E(\mathcal{T})$}, where
$\mathcal{T}=\{T_i\mid i\in I\}$.

Clearly, $\mathcal{A}'\equiv\mathcal{A}_P$ realizing $p_\infty(x)$
is not elementary embeddable into $\mathcal{A}_P$ and can not be
represented as a disjoint $P$-combination of
$\mathcal{A}'_i\equiv\mathcal{A}_i$, $i\in I$. At the same time,
there are $E$-combinations such that all
$\mathcal{A}'\equiv\mathcal{A}_E$ can be represented as
$E$-combinations of some $\mathcal{A}'_j\equiv\mathcal{A}_i$. We
call this representability of $\mathcal{A}'$ to be the {\em
$E$-representability}.

If there is $\mathcal{A}'\equiv\mathcal{A}_E$ which is not
$E$-representable, we have the $E'$-representability replacing $E$
by $E'$ such that $E'$ is obtained from $E$ adding equivalence
classes with models for all theories $T$, where $T$ is a theory of
a restriction $\mathcal{B}$ of a structure
$\mathcal{A}'\equiv\mathcal{A}_E$ to some $E$-class and
$\mathcal{B}$ is not elementary equivalent to the structures
$\mathcal{A}_i$. The resulting structure $\mathcal{A}_{E'}$ (with
the $E'$-representability) is a {\em
$e$-completion}\index{$e$-completion}, or a {\em
$e$-saturation}\index{$e$-saturation}, of $\mathcal{A}_{E}$. The
structure $\mathcal{A}_{E'}$ itself is called {\em
$e$-complete}\index{Structure!$e$-complete}, or {\em
$e$-saturated}\index{Structure!$e$-saturated}, or {\em
$e$-universal}\index{Structure!$e$-universal}, or {\em
$e$-largest}\index{Structure!$e$-largest}.

For a structure $\mathcal{A}_E$ the number of {\em
new}\index{Structure!new} structures with respect to the
structures $\mathcal{A}_i$, i.~e., of the structures $\mathcal{B}$
which are pairwise elementary non-equivalent and elementary
non-equivalent to the structures $\mathcal{A}_i$, is called the
{\em $e$-spectrum}\index{$e$-spectrum} of $\mathcal{A}_E$ and
denoted by $e$-${\rm Sp}(\mathcal{A}_E)$.\index{$e$-${\rm
Sp}(\mathcal{A}_E)$} The value ${\rm sup}\{e$-${\rm
Sp}(\mathcal{A}'))\mid\mathcal{A}'\equiv\mathcal{A}_E\}$ is called
the {\em $e$-spectrum}\index{$e$-spectrum} of the theory ${\rm
Th}(\mathcal{A}_E)$ and denoted by $e$-${\rm Sp}({\rm
Th}(\mathcal{A}_E))$.\index{$e$-${\rm Sp}({\rm
Th}(\mathcal{A}_E))$} If structures $\mathcal{A}_i$ represent
theories $T_i$ of a family $\mathcal{T}$, consisting of $T_i$,
$i\in I$, then the $e$-spectrum $e$-${\rm Sp}(\mathcal{A}_E)$ is
denoted by $e$-${\rm Sp}(\mathcal{T})$.

If $\mathcal{A}_E$ does not have $E$-classes $\mathcal{A}_i$,
which can be removed, with all $E$-classes
$\mathcal{A}_j\equiv\mathcal{A}_i$, preserving the theory ${\rm
Th}(\mathcal{A}_E)$, then $\mathcal{A}_E$ is called {\em
$e$-prime}\index{Structure!$e$-prime}, or {\em
$e$-minimal}\index{Structure!$e$-minimal}.

For a structure $\mathcal{A}'\equiv\mathcal{A}_E$ we denote by
${\rm TH}(\mathcal{A}')$ the set of all theories ${\rm
Th}(\mathcal{A}_i)$\index{${\rm Th}(\mathcal{A}_i)$} of
$E$-classes $\mathcal{A}_i$ in $\mathcal{A}'$.

By the definition, an $e$-minimal structure $\mathcal{A}'$
consists of $E$-classes with a minimal set ${\rm
TH}(\mathcal{A}')$. If ${\rm TH}(\mathcal{A}')$ is the least for
models of ${\rm Th}(\mathcal{A}')$ then $\mathcal{A}'$ is called
{\em $e$-least}.\index{Structure!$e$-least}

\medskip
{\bf Definition} \cite{cl}. Let $\overline{\mathcal{T}}_\Sigma$ be
the set of all complete elementary theories of a relational
language $\Sigma$. For a set
$\mathcal{T}\subset\overline{\mathcal{T}}_\Sigma$ we denote by
${\rm Cl}_E(\mathcal{T})$ the set of all theories ${\rm
Th}(\mathcal{A})$, where $\mathcal{A}$ is a structure of some
$E$-class in $\mathcal{A}'\equiv\mathcal{A}_E$,
$\mathcal{A}_E={\rm Comb}_E(\mathcal{A}_i)_{i\in I}$, ${\rm
Th}(\mathcal{A}_i)\in\mathcal{T}$. As usual, if $\mathcal{T}={\rm
Cl}_E(\mathcal{T})$ then $\mathcal{T}$ is said to be {\em
$E$-closed}.\index{Set!$E$-closed}

The operator ${\rm Cl}_E$ of $E$-closure can be naturally extended
to the classes $\mathcal{T}\subset\overline{\mathcal{T}}$, where
$\overline{\mathcal{T}}$ is the union of all
$\overline{\mathcal{T}}_\Sigma$ as follows: ${\rm
Cl}_E(\mathcal{T})$ is the union of all ${\rm
Cl}_E(\mathcal{T}_0)$ for subsets
$\mathcal{T}_0\subseteq\mathcal{T}$, where new language symbols
with respect to the theories in $\mathcal{T}_0$ are empty.

For a set $\mathcal{T}\subset\overline{\mathcal{T}}$ of theories
in a language $\Sigma$ and for a sentence $\varphi$ with
$\Sigma(\varphi)\subseteq\Sigma$ we denote by
$\mathcal{T}_\varphi$\index{$\mathcal{T}_\varphi$} the set
$\{T\in\mathcal{T}\mid\varphi\in T\}$. Any set
$\mathcal{T}_\varphi$ is called the {\em $\varphi$-neighbourhood},
or simply a {\em neighbourhood}, for $\mathcal{T}$, or the
($\varphi$-){\em definable} subset of $\mathcal{T}$. The set
$\mathcal{T}_\varphi$ is also called ({\em formula-} or {\em
sentence-}){\em definable} (by the sentence $\varphi$) with
respect to $\mathcal{T}$, or ({\em sentence-})$\mathcal{T}$-{\em
definable}, or simply {\em $s$-definable}.

\medskip
{\bf Proposition 1.1} \cite{cl}. {\em If
$\mathcal{T}\subset\overline{\mathcal{T}}$ is an infinite set and
$T\in\overline{\mathcal{T}}\setminus\mathcal{T}$ then $T\in{\rm
Cl}_E(\mathcal{T})$ {\rm (}i.e., $T$ is an {\sl accumulation
point} for $\mathcal{T}$ with respect to $E$-closure ${\rm
Cl}_E${\rm )} if and only if for any formula $\varphi\in T$ the
set $\mathcal{T}_\varphi$ is infinite.}

\medskip
If $T$ is an accumulation point for $\mathcal{T}$ then we also say
that $T$ is an {\em accumulation point} for ${\rm
Cl}_E(\mathcal{T})$.

\medskip
{\bf Theorem 1.2} \cite{cl}. {\em For any sets
$\mathcal{T}_0,\mathcal{T}_1\subset\overline{\mathcal{T}}$, ${\rm
Cl}_E(\mathcal{T}_0\cup\mathcal{T}_1)={\rm
Cl}_E(\mathcal{T}_0)\cup{\rm Cl}_E(\mathcal{T}_1)$.}

\medskip
{\bf Definition} \cite{cl}. Let $\mathcal{T}_0$ be a closed set in
a topological space $(\mathcal{T},\mathcal{O}_E(\mathcal{T}))$,
where $\mathcal{O}_E(\mathcal{T})=\{\mathcal{T}\setminus{\rm
Cl}_E(\mathcal{T}')\mid\mathcal{T}'\subseteq\mathcal{T}\}$. A
subset $\mathcal{T}'_0\subseteq\mathcal{T}_0$ is said to be {\em
generating}\index{Set!generating} if $\mathcal{T}_0={\rm
Cl}_E(\mathcal{T}'_0)$. The generating set $\mathcal{T}'_0$ (for
$\mathcal{T}_0$) is {\em minimal}\index{Set!generating!minimal} if
$\mathcal{T}'_0$ does not contain proper generating subsets. A
minimal generating set $\mathcal{T}'_0$ is {\em
least}\index{Set!generating!least} if $\mathcal{T}'_0$ is
contained in each generating set for $\mathcal{T}_0$.

\medskip
{\bf Theorem 1.3} \cite{cl}. {\em If $\mathcal{T}'_0$ is a
generating set for a $E$-closed set $\mathcal{T}_0$ then the
following conditions are equivalent:

$(1)$ $\mathcal{T}'_0$ is the least generating set for
$\mathcal{T}_0$;

$(2)$ $\mathcal{T}'_0$ is a minimal generating set for
$\mathcal{T}_0$;

$(3)$ any theory in $\mathcal{T}'_0$ is isolated by some set
$(\mathcal{T}'_0)_\varphi$, i.e., for any $T\in\mathcal{T}'_0$
there is $\varphi\in T$ such that
$(\mathcal{T}'_0)_\varphi=\{T\}$;

$(4)$ any theory in $\mathcal{T}'_0$ is isolated by some set
$(\mathcal{T}_0)_\varphi$, i.e., for any $T\in\mathcal{T}'_0$
there is $\varphi\in T$ such that
$(\mathcal{T}_0)_\varphi=\{T\}$.}

\medskip

\medskip
{\bf Definition} \cite{at}. Let $\mathcal{T}$ be a family of
theories and $T$ be a theory, $T\notin\mathcal{T}$. The theory $T$
is called {\em $\mathcal{T}$-approximated}, or {\em approximated
by} $\mathcal{T}$, or {\em $\mathcal{T}$-approximable}, or a {\em
pseudo-$\mathcal{T}$-theory}, if for any formula $\varphi\in T$
there is $T'\in\mathcal{T}$ such that $\varphi\in T'$.

If $T$ is $\mathcal{T}$-approximated then $\mathcal{T}$ is called
an {\em approximating family} for $T$, theories $T'\in\mathcal{T}$
are {\em approximations} for $T$, and $T$ is an {\em accumulation
point} for $\mathcal{T}$.

An approximating family $\mathcal{T}$ is called {\em $e$-minimal}
if for any sentence $\varphi\in\Sigma(T)$, $\mathcal{T}_\varphi$
is finite or $\mathcal{T}_{\neg\varphi}$ is finite.

\medskip
It was shown in \cite{at} that any $e$-minimal family
$\mathcal{T}$ has unique accumulation point $T$ with respect to
neighbourhoods $\mathcal{T}_\varphi$, and $\mathcal{T}\cup\{T\}$
is also called {\em $e$-minimal}.

\medskip
Following \cite{RSrank} we define the {\em rank} ${\rm RS}(\cdot)$
for the families of theories, similar to Morley rank
\cite{Morley}, and a hierarchy with respect to these ranks in the
following way.

For the empty family $\mathcal{T}$ we put the rank ${\rm
RS}(\mathcal{T})=-1$, for finite nonempty families $\mathcal{T}$
we put ${\rm RS}(\mathcal{T})=0$, and for infinite families
$\mathcal{T}$~--- ${\rm RS}(\mathcal{T})\geq 1$.

For a family $\mathcal{T}$ and an ordinal $\alpha=\beta+1$ we put
${\rm RS}(\mathcal{T})\geq\alpha$ if there are pairwise
inconsistent $\Sigma(\mathcal{T})$-sentences $\varphi_n$,
$n\in\omega$, such that ${\rm
RS}(\mathcal{T}_{\varphi_n})\geq\beta$, $n\in\omega$.

If $\alpha$ is a limit ordinal then ${\rm
RS}(\mathcal{T})\geq\alpha$ if ${\rm RS}(\mathcal{T})\geq\beta$
for any $\beta<\alpha$.

We set ${\rm RS}(\mathcal{T})=\alpha$ if ${\rm
RS}(\mathcal{T})\geq\alpha$ and ${\rm
RS}(\mathcal{T})\not\geq\alpha+1$.

If ${\rm RS}(\mathcal{T})\geq\alpha$ for any $\alpha$, we put
${\rm RS}(\mathcal{T})=\infty$.

A family $\mathcal{T}$ is called {\em $e$-totally transcendental},
or {\em totally transcendental}, if ${\rm RS}(\mathcal{T})$ is an
ordinal.

Similarly \cite{Morley}, for a nonempty family $\mathcal{T}$, we
denote by $\mathcal{B}(\mathcal{T})$ the Boolean algebra
consisting of all subfamilies $\mathcal{T}_{\varphi}$, where
$\varphi$ are sentences in the language $\Sigma(\mathcal{T})$.

\medskip
{\bf Theorem 1.4} \cite{RSrank, Morley}. {\em A nonempty family
$\mathcal{T}$ is $e$-totally transcendental if and only if the
Boolean algebra $\mathcal{B}(\mathcal{T})$ is superatomic.}

\medskip
{\bf Proposition 1.5} \cite{RSrank}. {\em If an infinite family
$\mathcal{T}$ does not have $e$-minimal subfamilies
$\mathcal{T}_\varphi$ then $\mathcal{T}$ is not $e$-totally
transcendental.}

\medskip
If $\mathcal{T}$ is $e$-totally transcendental, with ${\rm
RS}(\mathcal{T})=\alpha\geq 0$, we define the {\em degree} ${\rm
ds}(\mathcal{T})$ of $\mathcal{T}$ as the maximal number of
pairwise inconsistent sentences $\varphi_i$ such that ${\rm
RS}(\mathcal{T}_{\varphi_i})=\alpha$.

\medskip
{\bf Proposition 1.6} \cite{RSrank}. {\em A family $\mathcal{T}$
is $e$-minimal if and only if ${\rm RS}(\mathcal{T})=1$ and ${\rm
ds}(\mathcal{T})=1$.}

\medskip
{\bf Proposition 1.7} \cite{RSrank}. {\em For any family
$\mathcal{T}$, ${\rm RS}(\mathcal{T})={\rm RS}({\rm
Cl}_E(\mathcal{T}))$, and if $\mathcal{T}$ is nonempty and
$e$-totally transcendental then ${\rm ds}(\mathcal{T})={\rm
ds}({\rm Cl}_E(\mathcal{T}))$.}

\medskip
Recall the definition of the Cantor--Bendixson rank. It is defined
on the elements of a topological space $X$ by induction: ${\rm
CB}_X(p)\geq 0$ for all $p\in X$; ${\rm CB}_X(p)\geq\alpha$ if and
only if for any $\beta<\alpha$, $p$ is an accumulation point of
the points of ${\rm CB}_X$-rank at least $\beta$. ${\rm
CB}_X(p)=\alpha$ if and only if both ${\rm CB}_X(p)\geq\alpha$ and
${\rm CB}_X(p)\ngeq\alpha+1$ hold; if such an ordinal $\alpha$
does not exist then ${\rm CB}_X(p)=\infty$. Isolated points of $X$
are precisely those having rank $0$, points of rank $1$ are those
which are isolated in the subspace of all non-isolated points, and
so on. For a non-empty $C\subseteq X$  we define ${\rm
CB}_X(C)=\sup\{{\rm CB}_X(p)\mid p\in C\}$; in this way ${\rm
CB}_X(X)$ is defined and ${\rm CB}_X(\{p\})={\rm CB}_X(p)$ holds.
If $X$ is compact and $C$ is closed in $X$ then the sup is
achieved: ${\rm CB}_X(C)$ is the maximum value of ${\rm CB}_X(p)$
for $p\in C$; there are  finitely many points of maximum rank in
$C$ and the number of such points is the \emph{${\rm
CB}_X$-degree} of $C$, denoted by $n_X(C)$.

If $X$ is countable and compact then ${\rm CB}_X(X)$ is a
countable ordinal and every  closed subset has ordinal-valued rank
and finite ${\rm CB}_X$-degree $n_X(X)\in\omega\setminus\{0\}$.

For any ordinal $\alpha$ the set $\{p\in X\mid{\rm
CB}_X(p)\geq\alpha\}$ is called the $\alpha$-th {\em ${\rm
CB}$-derivative} $X_\alpha$ of $X$.

Elements $p\in X$ with ${\rm CB}_X(p)=\infty$ form the {\em
perfect kernel} $X_\infty$ of $X$.

Clearly, $X_\alpha\supseteq X_{\alpha+1}$, $\alpha\in{\rm Ord}$,
and $X_\infty=\bigcap\limits_{\alpha\in{\rm Ord}}X_\alpha$.

It is noticed in \cite{RSrank} that any $e$-totally transcendental
family $\mathcal{T}$ defines a superatomic Boolean algebra
$\mathcal{B}(\mathcal{T})$ with ${\rm RS}(\mathcal{T})={\rm
CB}_{\mathcal{B}(\mathcal{T})}(B(\mathcal{T}))$, ${\rm
ds}(\mathcal{T})=n_{\mathcal{B}(\mathcal{T})}(B(\mathcal{T}))$,
i.e., the pair $({\rm RS}(\mathcal{T}),{\rm ds}(\mathcal{T}))$
consists of Cantor--Bendixson invariants for
$\mathcal{B}(\mathcal{T})$ \cite{BA}.

By the definition for any $e$-totally transcendental family
$\mathcal{T}$ each theory $T\in\mathcal{T}$ obtains the ${\rm
CB}$-rank ${\rm CB}_\mathcal{T}(T)$ starting with
$\mathcal{T}$-isolated points $T_0$, of ${\rm
CB}_\mathcal{T}(T_0)=0$. We will denote the values ${\rm
CB}_\mathcal{T}(T)$ by ${\rm RS}_\mathcal{T}(T)$ as the rank for
the point $T$ in the topological space on $\mathcal{T}$ which is
defined with respect to $\Sigma(\mathcal{T})$-sentences.

\medskip
{\bf Definition} \cite{RSrank}. Let $\alpha$ be an ordinal. A
family $\mathcal{T}$ of rank $\alpha$ is called {\em
$\alpha$-minimal} if for any sentence $\varphi\in\Sigma(T)$, ${\rm
RS}(\mathcal{T}_\varphi)<\alpha$ or ${\rm
RS}(\mathcal{T}_{\neg\varphi})<\alpha$.

\medskip
{\bf Proposition 1.8} \cite{RSrank}. {\em $(1)$ A family
$\mathcal{T}$ is $0$-minimal if and only if $\mathcal{T}$ is a
singleton.

$(2)$ A family $\mathcal{T}$ is $1$-minimal if and only if
$\mathcal{T}$ is $e$-minimal.

$(3)$ For any ordinal $\alpha$ a family $\mathcal{T}$ is
$\alpha$-minimal if and only if ${\rm RS}(\mathcal{T})=\alpha$ and
${\rm ds}(\mathcal{T})=1$. }

\medskip
{\bf Proposition 1.9} \cite{RSrank}. {\em For any family
$\mathcal{T}$, ${\rm RS}(\mathcal{T})=\alpha$, with ${\rm
ds}(\mathcal{T})=n$, if and only if $\mathcal{T}$ is represented
as a disjoint union of subfamilies
$\mathcal{T}_{\varphi_1},\ldots,\mathcal{T}_{\varphi_n}$, for some
pairwise inconsistent sentences $\varphi_1,\ldots,\varphi_n$, such
that each $\mathcal{T}_{\varphi_i}$ is $\alpha$-minimal.}

\section{Calculi for families of theories. Links for sentence-definable and diagram-definable families}

In this section we define calculi for families of theories,
similar to first-order calculi for sentences, as well as discuss
properties and links for these calculi.

For a family $\mathcal{T}$ and sentences $\varphi$ and $\psi$ we
say that $\varphi$ {\em $\mathcal{T}$-forces} $\psi$, written
$\varphi\vdash_\mathcal{T}\psi$ if $\mathcal{T}_\varphi\subseteq
\mathcal{T}_\psi$.

We put $\vdash_\mathcal{T}\psi$ if $\mathcal{T}_\psi=\mathcal{T}$,
and $\varphi\vdash_\mathcal{T}$ if
$\mathcal{T}_\varphi=\emptyset$. For $\vdash_\mathcal{T}\psi$ we
say that $\psi$ is {\em $\mathcal{T}$-provable}, and if
$\varphi\vdash_\mathcal{T}$ then we say that $\varphi$ is {\em
$\mathcal{T}$-contradictory} or {\em $\mathcal{T}$-inconsistent}.

By the definition the relation $\vdash_\mathcal{T}\psi$ is
equivalent to $\chi\vdash_\mathcal{T}\psi$ for any identically
true sentence $\chi$, and $\varphi\vdash_\mathcal{T}$ is
equivalent to $\varphi\vdash_\mathcal{T}\theta$ for any
identically false sentence $\theta$. So below we consider only
relations of form $\varphi\vdash_\mathcal{T}\psi$ and their
natural modifications.

Ordinary axioms and rules for calculi of sentences can be
naturally transformed for the relations
$\varphi\vdash_\mathcal{T}\psi$ obtaining {\em
$\mathcal{T}$-calculi}, i.e., calculi with respect to families
$\mathcal{T}$.

Clearly, $\varphi\vdash_\emptyset\psi$ for any sentences $\varphi$
and $\psi$. Therefore there are sentences $\varphi$ and $\psi$
such that $\varphi\vdash_\mathcal{T}\psi$ but
$\varphi\not\vdash\psi$. Indeed, if $\varphi$ and $\psi$ are
sentences in a language $\Sigma$ satisfying $\vdash\varphi$ and
$\not\vdash\psi$ then we have $\varphi\not\vdash\psi$ whereas
$\varphi\vdash_\emptyset\psi$. Besides, for the set
$\mathcal{T}_\Sigma$ of all theories in the language $\Sigma$ and
for $\mathcal{T}=(\mathcal{T}_\Sigma)_\psi$ we have
$\varphi\vdash_\mathcal{T}\psi$. Additionally, for any sentence
$\varphi$ which does not belong to theories in a family
$\mathcal{T}$, i.e., $\mathcal{T}_\varphi=\emptyset$, and for any
sentence $\psi$ we have $\varphi\vdash_\mathcal{T}\psi$.

The following obvious proposition asserts that the relation
$\varphi\vdash_\mathcal{T}\psi$ is monotone under $\vdash$ and
inclusion:

\medskip
{\bf Proposition 2.1.} {\em For any sentences
$\varphi,\varphi',\psi,\psi'$ and families
$\mathcal{T},\mathcal{T}'$, if $\varphi'\vdash\varphi$,
$\psi\vdash\psi'$, and $\mathcal{T}'\subseteq\mathcal{T}$ then
$\varphi\vdash_\mathcal{T}\psi$ implies
$\varphi'\vdash_{\mathcal{T}'}\psi'$.}

\medskip
The following proposition asserts the finite character for the
relations $\varphi\vdash_\mathcal{T}\psi$.

\medskip
{\bf Proposition 2.2.} {\em For any sentences $\varphi$, $\psi$
and a family $\mathcal{T}$ of theories the following conditions
are equivalent:

$(1)$ $\varphi\vdash_\mathcal{T}\psi$;

$(2)$ $\varphi\vdash_{\mathcal{T}_0}\psi$ for any finite
$\mathcal{T}_0\subseteq\mathcal{T}$;

$(3)$ $\varphi\vdash_{\{T\}}\psi$ for any singleton
$\{T\}\subseteq\mathcal{T}$.}

\medskip
Proof. The implications $(1)\Rightarrow(2)$ and
$(2)\Rightarrow(3)$ hold by Proposition 2.1.

$(3)\Rightarrow(1)$. In view of $\varphi\vdash_\emptyset\psi$ it
suffices to show $\varphi\vdash_\mathcal{T}\psi$ for nonempty
$\mathcal{T}$ having $\varphi\vdash_{\{T\}}\psi$ for any singleton
$\{T\}\subseteq\mathcal{T}$. But if $T\in\mathcal{T}_\varphi$ then
$T\in\{T\}_\varphi$ and using $\varphi\vdash_{\{T\}}\psi$ we
obtain $T\in\{T\}_\psi$ implying $T\in\mathcal{T}_\psi$. Thus,
$\varphi\vdash_\mathcal{T}\psi$.~$\Box$

\medskip
{\bf Proposition 2.3.} {\em For any sentences $\varphi$ and $\psi$
in a language $\Sigma$ the following conditions are equivalent:

$(1)$ $\varphi\vdash\psi$;

$(2)$ $\varphi\vdash_{\mathcal{T}_\Sigma}\psi$;

$(3)$ $\varphi\vdash_{\mathcal{T}}\psi$ for any {\rm (}finite{\rm
)} family {\rm (}singleton{\rm )} $\mathcal{T}\subseteq
\mathcal{T}_\Sigma$;

$(4)$ $\varphi\vdash_{\mathcal{T}}\psi$ for any {\rm (}finite{\rm
)} family {\rm (}singleton{\rm )} $\mathcal{T}$;

$(5)$ $T\cup\{\varphi\}\vdash\psi$ for any
$T\in\mathcal{T}_\Sigma$.}

\medskip
Proof. $(4)\Rightarrow(3)$ and $(3)\Rightarrow(2)$ are obvious
using Proposition 2.2.

$(2)\Rightarrow(1)$. Assume on contrary that
$\varphi\vdash_{\mathcal{T}_\Sigma}\psi$ and
$\varphi\not\vdash\psi$. Then $\varphi\wedge\neg\psi$ is
consistent. Extending $\{\varphi\wedge\neg\psi\}$ till a complete
theory $T$ in the language $\Sigma$ we obtain
$T\in(\mathcal{T}_\Sigma)_\varphi$ and
$T\notin(\mathcal{T}_\Sigma)_\psi$ contradicting
$\varphi\vdash_{\mathcal{T}_\Sigma}\psi$.

$(1)\Rightarrow(4)$. If $\varphi\vdash\psi$ then for any theory
$T$ with $\varphi\in T$ we have $\psi\in T$, hence
$\mathcal{T}_\varphi\subseteq \mathcal{T}_\psi$ for any family
$\mathcal{T}$, i.e., $\varphi\vdash_{\mathcal{T}}\psi$.

$(3)\Leftrightarrow(5)$. $\varphi\vdash_{\{T\}}\psi$ means that
$\varphi\in T$ implies $\psi\in T$. So if $\varphi\in T$ then
$T\vdash\psi$ implying $T\cup\{\varphi\}\vdash\psi$. Otherwise if
$\varphi\notin T$ then $\neg\varphi\in T$. Therefore we have
$\{\varphi,\neg\varphi\}\vdash\psi$ implying
$T\cup\{\varphi\}\vdash\psi$. Conversely, assuming on contrary
$\varphi\not\vdash_{\{T\}}\psi$ we have $\varphi\in T$ and
$\psi\notin T$, so $\neg\psi\in T$. Hence
$T\cup\{\varphi\}\not\vdash\psi$ since $T$ is complete theory and
containing $\neg\psi$ it can not force $\psi$, i.e., $T$ can not
contain $\psi$.~$\Box$

\medskip
{\bf Definition.} If $\mathcal{T}$ is a family of theories and
$\Phi$ is a set of sentences, then we put
$\mathcal{T}_\Phi=\bigcap\limits_{\varphi\in\Phi}\mathcal{T}_\varphi$
and the set $\mathcal{T}_\Phi$ is called ({\em type-} or {\em
diagram-}){\em definable} (by the set $\Phi$) with respect to
$\mathcal{T}$, or ({\em diagram-})$\mathcal{T}$-{\em definable},
or simply {\em $d$-definable}.
\medskip

By the definition we have the following properties:

\medskip
0. Any $d$-definable subfamily of $E$-closed family $\mathcal{T}$
is again $E$-closed.

\medskip
1. $\mathcal{T}_{\{\varphi\}}=\mathcal{T}_\varphi$.

\medskip
2. $\mathcal{T}_\Phi=\{T\in\mathcal{T}\mid\Phi\subseteq T\}$.

\medskip
3. $\mathcal{T}_\Phi=\mathcal{T}$ if and only if
$\Phi\subseteq\cap\mathcal{T}$. In particular,
$\mathcal{T}_\emptyset=\mathcal{T}$.

\medskip
4.
$\mathcal{T}_{\Phi\cup\Psi}=\mathcal{T}_\Phi\cap\mathcal{T}_\Psi$.

\medskip
5. $\mathcal{T}_\Phi=(\mathcal{T}_\Phi)_\Psi$ for any $\Psi$
consisting of sentences $\psi$ with $\Phi\vdash\psi$. In
particular, the operation $(\cdot)_\Phi$ is idempotent:
$(\mathcal{T}_\Phi)_\Phi=\mathcal{T}_\Phi$.

\medskip
6.
$\mathcal{T}_{\{\varphi_1,\ldots,\varphi_n\}}=\mathcal{T}_{\varphi_1\wedge\ldots\wedge\varphi_n}$,
i.e., definable sets $\mathcal{T}_\Phi$ by finite $\Phi$ are
sentence-definable.

\medskip
7. $\mathcal{T}_\Phi=\mathcal{T}_\Psi$, where $\Psi$ is the
closure of $\Phi$ under conjunctions.

\medskip
By the latter property, studying $d$-definable sets, we will
usually consider sets $\Phi$ closed under conjunctions. Moreover,
by Property 5, considering $d$-definable families we can
additionally assume that any $\Phi$ is closed under logical
conclusions with respect to $\vdash$. It means that it suffices to
assume that $\Phi$ corresponds a filter with respect to the family
of $d$-definable subsets of $\mathcal{T}$.

\medskip
8. For any sets $\Phi$ and $\Psi$ containing all their logical
conclusions,
$\mathcal{T}_{\Phi\cap\Psi}=\mathcal{T}_\Phi\cup\mathcal{T}_\Psi$.

\medskip
Indeed, if $T\in\mathcal{T}_{\Phi\cap\Psi}$ then
$\Phi\cap\Psi\subseteq T$. Assuming
$T\notin\mathcal{T}_\Phi\cup\mathcal{T}_\Psi$ we have
$\Phi\not\subseteq T$ and $\Psi\not\subseteq T$. So there are
sentences $\varphi\in\Phi\setminus T$ and $\psi\in\Psi\setminus
T$. Then $\varphi\vee\psi\notin T$. But by conjecture,
$\varphi\vee\psi\in\Phi\cap\Psi$ contradicting
$\Phi\cap\Psi\subseteq T$. Conversely, if
$T\in\mathcal{T}_\Phi\cup\mathcal{T}_\Psi$ then $\Phi\subseteq T$
or $\Psi\subseteq T$ implying $\Phi\cap\Psi\subseteq T$ and
$T\in\mathcal{T}_{\Phi\cap\Psi}$.

\medskip
9. For any $T\in\mathcal{T}$ and $\Phi\subseteq T$ with
$\Phi\vdash \varphi$ for all $\varphi\in T$,
$\mathcal{T}_\Phi=\{T\}$. So any set of axioms for $T$ isolates
$T$ in $\mathcal{T}$. In particular, since $T$ is an ultrafilter
and axiomatized by itself, $\mathcal{T}_T=\{T\}$.

\medskip
The following proposition gives obvious criteria for $d$-definable
sets to be $s$-definable.

\medskip
{\bf Proposition 2.4.} {\em For any $d$-definable set
$\mathcal{T}=\mathcal{T}_\Phi$, where $\Phi$ is closed under
conjunctions, and a sentence $\varphi\in\Phi$ the following
conditions are equivalent:

$(1)$ $\mathcal{T}$ is $s$-definable by $\varphi$:
$\mathcal{T}=\mathcal{T}_\varphi$;

$(2)$ $\varphi\vdash_\mathcal{T}\psi$ for any $\psi\in\Phi$;

$(3)$ each $\psi\in\Phi$ with $\psi\vdash\varphi$ satisfies
$\mathcal{T}_\varphi= \mathcal{T}_\psi$;

$(3)$ there are no $T\in\mathcal{T}$ containing
$\varphi\wedge\neg\psi$ for any $\psi\in\Phi$.}

\medskip
The sentence $\varphi$ with $\mathcal{T}_\varphi\ne\emptyset$ and
satisfying the conditions in Proposition 2.4 is called {\em
$\mathcal{T}$-isolating}, {\em $\mathcal{T}$-principal} or {\em
$\mathcal{T}$-complete} for $\Phi$, and $\Phi$ is called {\em
$\mathcal{T}$-isolated} or {\em $\mathcal{T}$-principal}.

By Proposition 2.3, $\mathcal{T}_\Sigma$-isolating sentences are
isolating for $\Phi$, in the ordinary sense. Besides, if $\Phi$ is
forced by some $\varphi\in\Phi$ then for any family $\mathcal{T}$,
$\Phi$ is $\mathcal{T}$-isolated, but not vice versa.

\medskip
Clearly, each $d$-definable set $\mathcal{T}_\Phi$ equals the set
$\mathcal{T}_\theta$, where $\theta=\bigwedge\Phi$ with possibly
infinite conjunction and $\mathcal{T}_\theta$ is the set of all
theories $T\in\mathcal{T}$ containing conjunctive members of
$\theta$.

By Property 8 finite unions of $d$-definable sets are again
$d$-definable. Considering infinite unions $\mathcal{T}'$ of
$d$-definable sets $\mathcal{T}_{\Phi_i}$, $i\in I$, we can
represent them by sets of formulas with infinite disjunctions
$\bigvee\limits_{i\in I}\varphi_i$, $\varphi_i\in\Phi_i$. We call
these unions $\mathcal{T}'$ as {\em $d_\infty$-definable} sets.

Now the definability for subfamilies of $\mathcal{T}$ can be
extended for infinite unions, intersections and their complements.
Notice that since all singletons $\{T\}\subseteq\mathcal{T}$ are
$d$-definable, each subfamily $\mathcal{T}'\subseteq\mathcal{T}$
is $d_\infty$-definable.

\medskip
The relations $\varphi\vdash_\mathcal{T}\psi$ can be naturally
spread to sets $\Phi$ and $\Psi$ of sentences producing relations
$\Phi\vdash_\mathcal{T}\Psi$ meaning
$\mathcal{T}_\Phi\subseteq\mathcal{T}_\Psi$.

By Proposition 2.1 the relations $\Phi\vdash_\mathcal{T}\Psi$ are
again monotone:

\medskip
{\bf Proposition 2.5.} {\em For any sets $\Phi,\Phi',\Psi,\Psi'$
of sentences and families $\mathcal{T}, \mathcal{T}'$, if
$\Phi'\vdash\Phi$, $\Psi\vdash\Psi'$, and
$\mathcal{T}'\subseteq\mathcal{T}$ then
$\Phi\vdash_\mathcal{T}\Psi$ implies
$\Phi'\vdash_{\mathcal{T}'}\Psi'$.}

\medskip
Proposition 2.2 implies the following:

\medskip
{\bf Proposition 2.6.} {\em For any sets $\Phi$ and $\Psi$ of
sentences and a family $\mathcal{T}$ of theories the following
conditions are equivalent:

$(1)$ $\Phi\vdash_\mathcal{T}\Psi$;

$(2)$ $\Phi\vdash_{\mathcal{T}_0}\Psi$ for any finite
$\mathcal{T}_0\subseteq\mathcal{T}$;

$(3)$ $\Phi\vdash_{\{T\}}\Psi$ for any singleton
$\{T\}\subseteq\mathcal{T}$.}

\medskip
Proposition 2.3 immediately implies

\medskip
{\bf Proposition 2.7.} {\em For any sets $\Phi$ and $\Psi$ of
sentences in a language $\Sigma$ the following conditions are
equivalent:

$(1)$ $\Phi\vdash\Psi$, i.e., each sentence in $\Psi$ is forced by
some conjunction of sentences in $\Phi$;

$(2)$ $\Phi\vdash_{\mathcal{T}_\Sigma}\Psi$;

$(3)$ $\Phi\vdash_{\mathcal{T}}\Psi$ for any {\rm (}finite{\rm )}
family {\rm (}singleton{\rm )} $\mathcal{T}\subseteq
\mathcal{T}_\Sigma$;

$(4)$ $\Phi\vdash_{\mathcal{T}}\Psi$ for any {\rm (}finite{\rm )}
family {\rm (}singleton{\rm )} $\mathcal{T}$.}

\medskip
Extending the list for criteria of $\Phi\vdash_\mathcal{T}\Psi$ we
have the following:

\medskip
{\bf Theorem 2.8.} {\em For any sets $\Phi$ and $\Psi$ of
sentences and a family $\mathcal{T}$ of theories the following
conditions are equivalent:

$(1)$ $\Phi\vdash_\mathcal{T}\Psi$;

$(2)$ $\Phi\vdash_{{\rm Cl}_E(\mathcal{T})}\Psi$.}

\medskip
Proof. Since $\mathcal{T}\subseteq{\rm Cl}_E(\mathcal{T})$ we have
$(2)\Rightarrow(1)$ by Proposition 2.5.

$(1)\Rightarrow(2)$. Assume that $\Phi\vdash_\mathcal{T}\Psi$. It
suffices to show that if $\varphi\in\Phi$, $\psi\in\Psi$ with
$\mathcal{T}_\varphi\subseteq\mathcal{T}_\psi$ then $({\rm
Cl}_E(\mathcal{T}))_\varphi\subseteq({\rm
Cl}_E(\mathcal{T}))_\psi$. Let $T\in({\rm
Cl}_E(\mathcal{T}))_\varphi$ By the hypothesis we can assume that
$T\in {\rm Cl}_E(\mathcal{T})\setminus\mathcal{T}$ and using
Proposition 1.1 we have infinite $\mathcal{T}_\chi$ for any
$\chi\in T$. Since $\varphi\in T$,
$(\mathcal{T}_\varphi)_\chi=\mathcal{T}_{\varphi\wedge\chi}$ are
also infinite for any $\chi\in T$ and therefore
$\mathcal{T}_\varphi\subseteq\mathcal{T}_\psi$ implies that all
$(\mathcal{T}_\psi)_\chi$ are infinite. Thus again by Proposition
1.1, $T\in{\rm Cl}_E(\mathcal{T}_\psi)=({\rm
Cl}_E(\mathcal{T}))_\psi$.~$\Box$

\medskip
Theorem 2.8 immediately implies the following:

\medskip
{\bf Corollary 2.9.} {\em For any sets $\Phi$ and $\Psi$ of
sentences, and families $\mathcal{T}$, $\mathcal{T}'$,
$\mathcal{T}''$ of theories such that $\mathcal{T}'$ generates
${\rm Cl}_E(\mathcal{T})$ and
$\mathcal{T}'\subseteq\mathcal{T}''\subseteq{\rm
Cl}_E(\mathcal{T})$, the following conditions are equivalent:

$(1)$ $\Phi\vdash_\mathcal{T}\Psi$;

$(2)$ $\Phi\vdash_{\mathcal{T}'}\Psi$;

$(3)$ $\Phi\vdash_{\mathcal{T}''}\Psi$.}

\medskip
{\bf Remark 2.10.} Notice that in general case Corollary 2.9 can
not be extended to families $\mathcal{T}''\not\subseteq{\rm
Cl}_E(\mathcal{T})$. Indeed, taking any theory $T\notin{\rm
Cl}_E(\mathcal{T})$ we have, by Proposition 1.1, a sentence
$\chi\in T$ such that $({\rm Cl}_E(\mathcal{T}))_\chi$ is finite.
Since $T\notin({\rm Cl}_E(\mathcal{T}))_\chi$ and $({\rm
Cl}_E(\mathcal{T}))_\chi$ is finite, there is a sentence
$\theta\in T$ such that $({\rm
Cl}_E(\mathcal{T}))_\theta=\emptyset$. Thus for any inconsistence
sentence $\varphi$ we have $\theta\vdash_{{\rm
Cl}_E(\mathcal{T})}\varphi$ whereas $\theta\not\vdash_{{\rm
Cl}_E(\mathcal{T})\cup\{T\}}\varphi$.~$\Box$

\medskip
The assertions above show that for any family $\mathcal{T}$ there
are {\em calculi}, connected with ordinary calculi for first-order
sentences \cite{ErPa}, both for the relations
$\varphi\vdash_\mathcal{T}\psi$ and $\Phi\vdash_\mathcal{T}\Psi$,
which satisfy monotone properties, are reflexive
($\Phi\vdash_\mathcal{T}\Phi$) and transitive (if
$\Phi\vdash_\mathcal{T}\Psi$ and $\Psi\vdash_\mathcal{T}{\rm X}$
then $\Phi\vdash_\mathcal{T}{\rm X}$).

\medskip
{\bf Definition.} Sets $\Phi$ and $\Psi$ of sentences are called
{\em $\mathcal{T}$-equivalent}, written
$\Phi\equiv_{\mathcal{T}}\Psi$, if $\Phi\vdash_\mathcal{T}\Psi$
and $\Psi\vdash_\mathcal{T}\Phi$, i.e.,
$\mathcal{T}_\Phi=\mathcal{T}_\Psi$.

Sentences $\varphi$ and $\psi$ are called {\em
$\mathcal{T}$-equivalent}, written
$\varphi\equiv_{\mathcal{T}}\psi$, if
$\{\varphi\}\equiv_{\mathcal{T}}\{\psi\}$.

\medskip
Clearly, the relations $\equiv_{\mathcal{T}}$ are equivalent
relations both for sentences and for sets of sentences.

Proposition 2.7 and Theorem 2.8 immediately implies the following:

\medskip
{\bf Proposition 2.11.} {\em For any sets $\Phi$ and $\Psi$ of
sentences in a language $\Sigma$ the following conditions are
equivalent:

$(1)$ $\Phi\vdash\Psi$ and $\Psi\vdash\Phi$, i.e., $\Phi$ and
$\Psi$ force each other;

$(2)$ $\Phi\equiv_{\mathcal{T}_\Sigma}\Psi$;

$(3)$ $\Phi\equiv_{\mathcal{T}}\Psi$ for any {\rm (}finite{\rm )}
family {\rm (}singleton{\rm )} $\mathcal{T}\subseteq
\mathcal{T}_\Sigma$;

$(4)$ $\Phi\equiv_{\mathcal{T}}\Psi$ for any {\rm (}finite{\rm )}
family {\rm (}singleton{\rm )} $\mathcal{T}$.}

\medskip
{\bf Corollary 2.12.} {\em For any sentences $\varphi$ and $\psi$
in a language $\Sigma$ the following conditions are equivalent:

$(1)$ $\varphi\vdash\psi$ and $\psi\vdash\varphi$;

$(2)$ $\varphi\equiv_{\mathcal{T}_\Sigma}\psi$;

$(3)$ $\varphi\equiv_{\mathcal{T}}\psi$ for any {\rm (}finite{\rm
)} family {\rm (}singleton{\rm )} $\mathcal{T}\subseteq
\mathcal{T}_\Sigma$;

$(4)$ $\varphi\equiv_{\mathcal{T}}\psi$ for any {\rm (}finite{\rm
)} family {\rm (}singleton{\rm )} $\mathcal{T}$.}

\medskip
Theorem 2.8 implies

\medskip
{\bf Corollary 2.13.} {\em For any sets $\Phi$ and $\Psi$ of
sentences, and families $\mathcal{T}$, $\mathcal{T}'$,
$\mathcal{T}''$ of theories such that $\mathcal{T}'$ generates
${\rm Cl}_E(\mathcal{T})$ and
$\mathcal{T}'\subseteq\mathcal{T}''\subseteq{\rm
Cl}_E(\mathcal{T})$, the following conditions are equivalent:

$(1)$ $\Phi\equiv_\mathcal{T}\Psi$;

$(2)$ $\Phi\equiv_{\mathcal{T}'}\Psi$;

$(3)$ $\Phi\equiv_{\mathcal{T}''}\Psi$.}

\section{Compactness and $E$-closed families}

{\bf Definition.} A $d$-definable set $\mathcal{T}_\Phi$ is called
{\em $\mathcal{T}$-consistent} if $\mathcal{T}_\Phi\ne\emptyset$,
and $\mathcal{T}_\Phi$ is called {\em locally
$\mathcal{T}$-consistent} if for any finite $\Phi_0\subseteq\Phi$,
$\mathcal{T}_{\Phi_0}$ is $\mathcal{T}$-consistent.

\medskip
Clearly, locally $\mathcal{T}$-consistent $\mathcal{T}$-principal
sets $\mathcal{T}_\Phi$ are $\mathcal{T}$-consistent.

Notice also that there are locally $\mathcal{T}$-consistent
$d$-definable sets $\mathcal{T}_\Phi$ which are not
$\mathcal{T}$-consistent. Indeed, let, for instance, $\mathcal{T}$
be an $e$-minimal family which does not contain its unique
accumulation point $T$. Then by the definition of accumulation
point, $\mathcal{T}_T$ is locally $\mathcal{T}$-consistent whereas
$\mathcal{T}_T=\emptyset$.

The following {\em Compactness Theorem} shows that this effect
does not occur for $E$-closed families.

\medskip
{\bf Theorem 3.1.} {\em For any nonempty $E$-closed family
$\mathcal{T}$, every locally $\mathcal{T}$-consistent
$d$-definable set $\mathcal{T}_\Phi$ is $\mathcal{T}$-consistent.}

\medskip
Proof. If all neighbourhoods $\mathcal{T}_\varphi$,
$\varphi\in\Phi$, contain same theory $T\in\mathcal{T}$ then
$\mathcal{T}_\Phi$ is $\mathcal{T}$-consistent. So we can assume
that $\Phi$ is infinite, closed under conjunctions,
non-$\mathcal{T}$-principal, and for any $\varphi\in\Phi$,
$\mathcal{T}_{\varphi}$ contains infinitely many theories in
$\mathcal{T}$. Now we extend step-by-step the set $\Phi$ till a
non-principal ultrafilter $T$ of sentences of the language
$\Sigma(\mathcal{T})$ such that each $\psi\in T$ satisfies
$|\mathcal{T}_{\psi}|\geq\omega$. Applying Proposition 1.1 we
obtain $T\in{\rm Cl}_E(\mathcal{T})=\mathcal{T}$, and by
$T\supset\Phi$ we have $T\in\mathcal{T}_\Phi$, i.e.,
$\mathcal{T}_\Phi$ is $\mathcal{T}$-consistent.~$\Box$

\medskip
Theories $T\in\mathcal{T}$ belonging to locally
$\mathcal{T}$-consistent $d$-definable sets $\mathcal{T}_\Phi$ are
called their {\em realizations}.

The following proposition, along Proposition 1.1 and compactness
above, clarifies the mechanism of construction of ${\rm
Cl}_E(\mathcal{T})$ via realizations of $d$-definable subfamilies
of $\mathcal{T}$.

\medskip
{\bf Proposition 3.2.} {\em For any family $\mathcal{T}$, ${\rm
Cl}_E(\mathcal{T})$ consists of elements of $\mathcal{T}$ and of
accumulation points realizing locally $\mathcal{T}$-consistent
$d$-definable sets $\mathcal{T}_\Phi$.}

\medskip
Proof. By monotonicity property in Proposition 2.5, implying
$\mathcal{T}_\Phi\supseteq\mathcal{T}_\Psi$ for
$\Phi\subseteq\Psi$, it suffices to note that for any theory $T$,
$T\in{\rm Cl}_E(\mathcal{T})$ if and only if $T$ is a (unique)
realization of locally $\mathcal{T}$-consistent $d$-definable
subfamily $\mathcal{T}_T$.~$\Box$

\medskip
The following theorem gives a criterion of existence of
$d$-definable family which is not $s$-definable.

\medskip
{\bf Theorem 3.3.} {\em For any $E$-closed family $\mathcal{T}$,
there is a $d$-definable family $\mathcal{T}_\Phi$ which is not
$s$-definable if and only if $\mathcal{T}$ is infinite.}

\medskip
Proof. If $\mathcal{T}$ is finite then each theory
$T\in\mathcal{T}$ is isolated by some sentence $\varphi$. So each
nonempty subfamily of $\mathcal{T}$ is $s$-definable by some
disjunction of the sentences $\varphi$. Thus, since the empty
subfamily of $\mathcal{T}$ is $s$-definable, by an inconsistent
sentence, then each $d$-definable family $\mathcal{T}_\Phi$ is
$s$-definable.

Now we assume that $\mathcal{T}$ is infinite. By compactness,
since $\mathcal{T}$ is $E$-closed and infinite, the set $\Phi$ of
all sentences $\varphi$ such that $|\mathcal{T}_{\neg\varphi}|=1$
is $\mathcal{T}$-consistent. Taking an arbitrary theory
$T\in\mathcal{T}_\Phi$ we obtain a $d$-definable singleton
$\mathcal{T}_T=\{T\}$ which can not be $s$-definable by choice of
$\Phi$.~$\Box$

\medskip
{\bf Remark 3.4.} Theorem 3.3 does not hold for families
$\mathcal{T}$ which are not $E$-closed. Indeed, take an arbitrary
$e$-minimal family $\mathcal{T}$, which does not contain its
(unique) accumulation point $T$. Repeating arguments for the proof
of Theorem 3.3 we find the set $\Phi$ which is locally
$\mathcal{T}$-consistent but $\mathcal{T}_\Phi=\emptyset$ in view
of $T\notin\mathcal{T}$. Since all $s$-definable subfamilies of
$\mathcal{T}$ are either finite or cofinite, the only possibility
for new $d$-definable subfamily of $\mathcal{T}$ is
$\mathcal{T}_\Phi$. Since $\mathcal{T}_\Phi$ is empty,
$\mathcal{T}$ does not have $d$-definable subfamilies which are
not $s$-definable.~$\Box$

\section{Dynamics of ranks with respect to definable subfamilies of theories}

Let $\mathcal{T}$ be a family of theories, $\Phi$ be a set of
sentences, $\alpha$ be an ordinal $\leq{\rm RS}(\mathcal{T})$ or
$-1$. The set $\Phi$ is called {\em $\alpha$-ranking} for
$\mathcal{T}$ if ${\rm RS}(\mathcal{T}_\Phi)=\alpha$. A sentence
$\varphi$ is called {\em $\alpha$-ranking} for $\mathcal{T}$ if
${\rm RS}(\mathcal{T}_{\{\varphi\}})=\alpha$.

The set $\Phi$ (the sentence $\varphi$) is called {\em ranking}
for $\mathcal{T}$ if it is $\alpha$-ranking for $\mathcal{T}$ with
some $\alpha$.

\medskip
{\bf Definition} \cite{at}. For a family $\mathcal{T}$, a theory
$T$ is {\em $\mathcal{T}$-finitely axiomatizable}, or {\em
finitely axiomatizable with respect to $\mathcal{T}$}, or {\em
$\mathcal{T}$-relatively finitely axiomatizable}, if
$\mathcal{T}_\varphi=\{T\}$ for some
$\Sigma(\mathcal{T})$-sentence $\varphi$.

For a family $\mathcal{T}$ of a language $\Sigma$, a sentence
$\varphi$ of the language $\Sigma$ is called {\em
$\mathcal{T}$-complete} if $\varphi$ isolates a unique theory in
$\mathcal{T}$, i.e., $\mathcal{T}_\varphi$ is a singleton.

\medskip
{\bf Proposition 4.1.} $(1)$ {\em A set $\Phi$ {\rm (}a sentence
$\varphi${\rm )} is $(-1)$-ranking for $\mathcal{T}$ if and only
if $\mathcal{T}=\emptyset$ or $\Phi$ {\rm (}respectively
$\varphi${\rm )} is inconsistent with theories in $\mathcal{T}$.

$(2)$ A set $\Phi$ {\rm (}a sentence $\varphi${\rm )} is
$0$-ranking for $\mathcal{T}$, with ${\rm
ds}(\mathcal{T}_\Phi)=n$, if and only if $\Phi$ {\rm
(}respectively $\varphi${\rm )} is consistent exactly with some
$n\in\omega\setminus\{0\}$ theories in $\mathcal{T}$.

$(3)$ Any $0$-ranking sentence $\varphi$ for $\mathcal{T}$, with
${\rm ds}(\mathcal{T}_\varphi)=n$, is $\mathcal{T}$-equivalent to
a disjunction of $n$ {\rm (}pairwise inconsistent{\rm )}
$\mathcal{T}$-complete sentences.}

\medskip
Proof. (1) and (2) immediately follow from the definition.

(3) In view of ${\rm RS}(\mathcal{T}_\varphi)=0$ and ${\rm
ds}(\mathcal{T}_\varphi)=n$ we have
$\mathcal{T}_\varphi=\{T_1,\ldots,T_n\}$ for some distinct
theories $T_1,\ldots,T_n\in\mathcal{T}$. Since the theories $T_i$
are distinct, there are sentences $\psi_i\in T_i$ such that
$\neg\psi_i\in T_j$ for $j\ne i$. Thus the formulas
$$(\varphi\wedge\psi_1\wedge\neg\psi_2\wedge\ldots\wedge\neg\psi_n),$$
$$\ldots,$$
$$(\varphi\wedge\neg\psi_1\wedge\neg\psi_2\wedge\ldots\wedge\neg\psi_{n-2}\wedge\psi_{n-1}\wedge\neg\psi_{n}),$$
$$(\varphi\wedge\neg\psi_1\wedge\ldots\wedge\neg\psi_{n-1})$$ are
$\mathcal{T}$-complete, pairwise inconsistent and such that their
disjunction is $\mathcal{T}$-equivalent to $\varphi$.~$\Box$

\medskip
{\bf Remark 4.2.} By Proposition 4.1, if $T\in\mathcal{T}$ then
$\Phi=T$ is $0$-ranking, with $\mathcal{T}_T=\{T\}$. More
generally, for any distinct $T_1,\ldots,T_n\in\mathcal{T}$ the set
$T_1\vee\ldots\vee
T_n=\{\varphi_1\vee\ldots\vee\varphi_n\mid\varphi_i\in T_i\}$ is
$0$-ranking, with ${\rm ds}(\mathcal{T}_{T_1\vee\ldots\vee
T_n})=n$.

\medskip
As shown in Remark 4.2 each finite subset
$\mathcal{T}_0\subseteq\mathcal{T}$ is $d$-definable, and
Proposition 4.1 gives a characterization for $\mathcal{T}_0$ to be
$s$-definable.

The following theorem produces a characterization for a subfamily
$\mathcal{T}'\subseteq\mathcal{T}$ to be $d$-definable.

\medskip
{\bf Theorem 4.3.} {\em A subfamily
$\mathcal{T}'\subseteq\mathcal{T}$ is $d$-definable in
$\mathcal{T}$ if and only if $\mathcal{T}'$ is $E$-closed in
$\mathcal{T}$, i.e., $\mathcal{T}'={\rm Cl}_E(\mathcal{T}')\cap
\mathcal{T}$.}

\medskip
Proof. In view of Remark 4.2 we can assume that $\mathcal{T}'$ is
infinite. Let $\mathcal{T}'$ be $d$-definable, i.e.,
$\mathcal{T}'=\mathcal{T}_\Phi$ for some set $\Phi$. By
Proposition 1.1, all theories in ${\rm Cl}_E(\mathcal{T}')$
contain the set $\Phi$, i.e., ${\rm Cl}_E(\mathcal{T}')\cap
\mathcal{T}\subseteq\mathcal{T}_\Phi$. Indeed, if a theory
$T\in{\rm Cl}_E(\mathcal{T}')$ does not contain a sentence
$\varphi\in\Phi$ then $\neg\varphi\in T$ and
$(\mathcal{T}')_{\neg\varphi}=\emptyset$ contradicting $T\in{\rm
Cl}_E(\mathcal{T}')$. Since $\mathcal{T}'\subseteq{\rm
Cl}_E(\mathcal{T}')\cap \mathcal{T}$, we have $\mathcal{T}'={\rm
Cl}_E(\mathcal{T}')\cap \mathcal{T}$, i.e., the subfamily
$\mathcal{T}'$ is $E$-closed in $\mathcal{T}$.

Now let the subfamily $\mathcal{T}'$ be $E$-closed in
$\mathcal{T}$. Denote by $\Phi$ the set $\bigcap\mathcal{T}'$,
i.e., the set of all $\Sigma(\mathcal{T})$-sentences belonging to
all theories in $\mathcal{T}'$. Clearly,
$\mathcal{T}'\subseteq\mathcal{T}_\Phi$. If
$\mathcal{T}'\subset\mathcal{T}_\Phi$, i.e., there is
$T\in\mathcal{T}_\Phi\setminus\mathcal{T}'$ then $T\notin{\rm
Cl}_E(\mathcal{T}')$. Applying Proposition 1.1 we find a sentence
$\varphi\in T$ such that $(\mathcal{T}')_\varphi$ is finite, say,
$(\mathcal{T}')_\varphi=\{T_1,\ldots,T_n\}$. Since $T_i\ne T$
there are sentences $\psi_i\in T\setminus T_i$, $i=1,\ldots,n$.
For the sentence
$\chi=\varphi\wedge\psi_1\wedge\ldots\wedge\psi_n$ we have
$\chi\in T$ and $(\mathcal{T}')_\chi=\emptyset$. It implies
$\neg\chi\in\Phi$, contradicting $T\in\mathcal{T}_\Phi$.~$\Box$

\medskip
The following proposition shows that $s$-definable subsets of a
family $\mathcal{T}$ witnessing ${\rm RS}(\mathcal{T})=\beta$
produce a hierarchy of $\alpha$-ranking sentences for all ordinals
$\alpha\leq\beta$.

\medskip
{\bf Proposition 4.4.} {\em For any ordinals $\alpha\leq\beta$, if
${\rm RS}(\mathcal{T})=\beta$ then ${\rm
RS}(\mathcal{T}_\varphi)=\alpha$ for some {\rm
(}$\alpha$-ranking{\rm )} sentence $\varphi$. Moreover, there are
${\rm ds}(\mathcal{T})$ pairwise $\mathcal{T}$-inconsistent
$\beta$-ranking sentences for $\mathcal{T}$, and if $\alpha<\beta$
then there are infinitely many pairwise $\mathcal{T}$-inconsistent
$\alpha$-ranking sentences for $\mathcal{T}$.}

\medskip
Proof. Since the  Boolean algebra $F(\mathcal{T})$ is superatomic
by Theorem 1.1, each $\mathcal{T}_\varphi$ belongs to a hierarchy
with respect to the rank ${\rm RS}(\cdot)$ starting with
singletons, $e$-minimal subfamilies, etc. Thus, each
$\mathcal{T}_\varphi$ obtains a value ${\rm
RS}(\mathcal{T}_\varphi)=\alpha$ in this hierarchy such that all
$\alpha\leq\beta$ are witnessed by some $\mathcal{T}_\varphi$. By
the definition of ${\rm RS}(\cdot)$, $\mathcal{T}$ can be divided
onto ${\rm ds}(\mathcal{T})$ disjoint parts $\mathcal{T}_\varphi$
having the rank $\beta$. Again by the definition, if
$\alpha<\beta$ then there are infinitely many pairwise
$\mathcal{T}$-inconsistent $\alpha$-ranking sentences for
$\mathcal{T}$.~$\Box$

\medskip
By Proposition 4.4, for every family $\mathcal{T}$ with ${\rm
RS}(\mathcal{T})=\beta\geq 0$ the possibilities for ${\rm
RS}(\mathcal{T}')$ with $\mathcal{T}'\subseteq \mathcal{T}$ are
realized by $s$-definable subsets $\mathcal{T}_\varphi$ with ${\rm
RS}(\mathcal{T}_\varphi)=\alpha$ for all $\alpha\leq\beta$. Thus
the following natural question arises for families $\mathcal{T}$
which are not $e$-totally transcendental.

\medskip
{\bf Question.} {\em Let $\mathcal{T}$ be a family with ${\rm
RS}(\mathcal{T})=\infty$. What are the ${\rm RS}$-possibilities
for $s$-definable / $d$-definable subfamilies of $\mathcal{T}$?}

\medskip
As shown in Remark 4.2 every finite subset of $\mathcal{T}$ is
$d$-definable and $0$-ranking. So in fact the question arises for
$\alpha\geq 0$ with $s$-definable subfamilies, and for $\alpha>0$
with $d$-definable subfamilies.

Partially answering the question we notice that obtaining an
$s$-definable / $d$-definable subfamily $\mathcal{T}_\beta$ of
$\mathcal{T}$ with ${\rm RS}(\mathcal{T}_\beta)=\beta\geq 0$ we
have, by Proposition 4.4, $s$-definable / $d$-definable
subfamilies $\mathcal{T}_\alpha$ of $\mathcal{T}$ with ${\rm
RS}(\mathcal{T}_\alpha)=\alpha$, for all ordinals
$\alpha\leq\beta$. Thus, the required ordinals $\alpha$ form an
initial segment.

\medskip
Illustrating the question we notice that, in some more or less
general cases, the possibility for $\alpha=0$ with $s$-definable
subfamilies can be realized:

\medskip
{\bf Remark 4.5.} If a family $\mathcal{T}$ has an
$\alpha$-ranking sentence, for $\alpha\geq 0$, it does not imply
that $\mathcal{T}$ is $e$-totally transcendental. Indeed, any
family $\mathcal{T}$, for instance, of functional language,
$e$-totally transcendental or not, and with a theory $T$ of an
one-element algebra has a $0$-ranking sentence $\varphi$ saying
that the universe is a singleton. Clearly,
$\mathcal{T}_\varphi=\{T\}$.

At the same time there are many examples of families of theories
without nonempty $s$-definable $e$-totally transcendental
subfamilies. Indeed, taking, for instance, a family
$\mathcal{T}_\Sigma$ of all theories in a language $\Sigma$
containing infinitely many predicate symbols, we can not control,
by a sentence, links between all predicates. In particular, there
are at least continuum many possibilities arbitrarily varying
empty/nonempty predicates. These variations produce unbounded
ranks for any nonempty $s$-definable subfamilies
$\mathcal{T}_\varphi$ implying ${\rm
RS}(\mathcal{T}_\varphi)=\infty$.

\medskip
The following theorem gives an answer to the question for
$d$-definable subfamilies of theories in countable languages. The
arguments for this answer can be naturally spread for arbitrary
languages.

\medskip
{\bf Theorem 4.6.} {\em Let $\mathcal{T}$ be a family of a
countable language $\Sigma$ and with ${\rm
RS}(\mathcal{T})=\infty$, $\alpha\in\{0,1\}$,
$n\in\omega\setminus\{0\}$. Then there is a $d$-definable
subfamily $\mathcal{T}_\Phi$ such that ${\rm
RS}(\mathcal{T}_\Phi)=\alpha$ and ${\rm ds}(\mathcal{T}_\Phi)=n$.}

\medskip
Proof. We fix a family $\mathcal{T}$ of countable language
$\Sigma$, with ${\rm RS}(\mathcal{T})=\infty$, a countable ordinal
$\alpha$, and $n\in\omega\setminus\{0\}$. By Theorem 4.3 it
suffices to find an $E$-closed subfamily $\mathcal{T}'$ in
$\mathcal{T}$ with ${\rm RS}(\mathcal{T}')=1$ and ${\rm
ds}(\mathcal{T}')=n$.

If $\alpha=0$ then $\mathcal{T}'$ exists by Remark 4.2.

If $\alpha=1$ we take $n$ pairwise inconsistent sentences
$\varphi_i$, $i=1,\ldots,n$, such that ${\rm
RS}(\mathcal{T}_{\varphi_i})=\infty$ and for each
$\mathcal{T}_{\varphi_i}$ find $E$-closed, in
$\mathcal{T}_{\varphi_i}$ (and so in $\mathcal{T}$), $e$-minimal
subfamily $\mathcal{T}'_i$ in the following way. We enumerate the
set of all $\Sigma$-sentences which force $\varphi_i$:
$\psi_{ik}$, $k\in\omega$, and form $\mathcal{T}'_i$ step-by-step
with respect to that enumeration using the following subfamilies
$\mathcal{T}_{ik}$ of $\mathcal{T}_{\varphi_i}$ with
$\mathcal{T}_{ik}\supseteq\mathcal{T}_{i,k+1}$.

At the initial step if $\mathcal{T}_{\psi_{i0}}$ is cofinite in
$\mathcal{T}_{\varphi_i}$, we set
$\mathcal{T}_{i0}=\mathcal{T}_{\psi_{i0}}$, if
$\mathcal{T}_{\varphi_i}=\mathcal{T}_{\psi_{i0}}$, and
$\mathcal{T}_{i0}=\mathcal{T}_{\psi_{i0}}\cup\{T_0\}$, with an
arbitrary theory
$T^i_0\in\mathcal{T}_{i}\setminus\mathcal{T}_{\psi_{i0}}$, if
$\mathcal{T}_{i}\ne\mathcal{T}_{\psi_{i0}}$. If
$\mathcal{T}_{\psi_{i0}}$ is co-infinite we repeat the process
replacing $\psi_{i0}$ by $\varphi_i\wedge\neg\psi_{i0}$: for
infinite $\mathcal{T}_{\varphi_i\wedge\neg\psi_{i0}}$ instead of
$\mathcal{T}_{\psi_{i0}}$.

Let at the step $k$ a family $\mathcal{T}_{ik}$ is already formed
with some theories $T^i_0,\ldots,T^i_r$ added to families
$\mathcal{T}_{\psi_{is}}$ or
$\mathcal{T}_{\varphi_i\wedge\neg\psi_{is}}$. Now we consider the
sentence $\psi_{i,k+1}$. If $(\mathcal{T}_{ik})_{\psi_{i,k+1}}$ is
cofinite in $\mathcal{T}_{ik}$, we set
$\mathcal{T}_{i,k+1}=\mathcal{T}_{ik}$, if
$\mathcal{T}_{ik}=(\mathcal{T}_{ik})_{\psi_{i,k+1}}$ modulo
$\{T^i_0,\ldots,T^i_r\}$, and
$\mathcal{T}_{i,k+1}=(\mathcal{T}_{ik})_{\psi_{i,k+1}}\cup\{T^i_{r+1}\}$,
with an arbitrary theory
$T^i_{r+1}\in\mathcal{T}_{i,k+1}\setminus((\mathcal{T}_{ik})_{\psi_{i,k+1}}\cup\{T^i_0,\ldots,T^i_r\})$,
if $\mathcal{T}_{ik}\ne(\mathcal{T}_{ik})_{\psi_{i,k+1}}$ modulo
$\{T^i_0,\ldots,T^i_r\}$. If $(\mathcal{T}_{ik})_{\psi_{i,k+1}}$
is co-infinite we repeat the process for infinite
$(\mathcal{T}_{ik})_{\varphi_i\wedge\neg\psi_{i,k+1}}$.

By the construction the subfamilies $\mathcal{T}'_i$ consisting of
the theories $T^i_0,\ldots,$ $T^i_r,$ $\ldots$, are infinite and
can not be divided into two infinite parts by $\Sigma$-sentences.
Indeed, $\mathcal{T}'_i$ is infinite because each set
$\{T^i_0,\ldots,T^i_r\}$ is extended in some step by some new
theory since
$\mathcal{T}_{ik}\ne(\mathcal{T}_{ik})_{\psi_{i,k+1}}$ modulo
$\{T^i_0,\ldots,T^i_r\}$ for some $\psi_{i,k+1}$ negating all
theories in $\{T^i_0,\ldots,T^i_r\}$ and some theory in
$\mathcal{T}_{ik}$. The subfamilies $\mathcal{T}'_i$ are
$e$-minimal since each $\Sigma$-sentence is equivalent to some
$\psi_{ik}$ modulo $\varphi_i$ and each $\psi_{ik}$ can divide
only $\mathcal{T}_{i0},\ldots,\mathcal{T}_{i,k-1}$ modulo
$\{T^i_0,\ldots,T^i_r\}$.

Thus, the subfamilies $\mathcal{T}'_i$ of $\mathcal{T}$ are
$e$-minimal, $i=1,\ldots,n$. By Proposition 1.6 we have ${\rm
RS}(\mathcal{T}'_i)=1$ and ${\rm ds}(\mathcal{T}'_i)=1$, and by
Proposition 1.7 we can assume that the families $\mathcal{T}'_i$
are $E$-closed in $\mathcal{T}$. Hence, for
$\mathcal{T}'=\mathcal{T}'_1\cup\ldots\cup\mathcal{T}'_n$, which
is $d$-definable by Theorem 4.3, we have ${\rm
RS}(\mathcal{T}')=1$ and ${\rm ds}(\mathcal{T}')=n$.~$\Box$

\medskip
{\bf Remark 4.7.} Notice that the arguments in the proof of
Theorem 4.6 do not work for $\alpha\geq 2$ since taking infinitely
many disjoint $s$-definable infinite subfamilies
$\mathcal{T}_{\varphi_i}$ we can not guarantee that
$\varphi_i\not\vdash_\mathcal{T}\psi_j$ for infinitely many
$\mathcal{T}$-disjoint sentences $\psi_j$. Thus constructing
$d$-definable $e$-minimal subfamilies $\mathcal{T}_i$ of
$\mathcal{T}_{\varphi_i}$ it is possible to obtain ${\rm
RS}\left(\bigcup\limits_i\mathcal{T}_i\right)\geq 3$, not ${\rm
RS}\left(\bigcup\limits_i\mathcal{T}_i\right)=2$.

At the same time, constructing countably many $d$-definable
subfamilies $\mathcal{T}_i$ of $\mathcal{T}_{\varphi_i}$,
$i\in\omega$, with pairwise inconsistent $\varphi_i$, we can
choose some infinite $I\subseteq\omega$, such that accumulation
points $T_i$ for $\mathcal{T}_i$, $i\in I$, form an $e$-minimal
family. Thus, possibly loosing the $d$-definability we obtain a
$d_\infty$-definable subfamily $\mathcal{T}'=\bigcup\limits_{i\in
I}\mathcal{T}_i$ with ${\rm RS}(\mathcal{T}')=2$ and ${\rm
ds}(\mathcal{T}')=1$. Taking some $n$ disjoint $\mathcal{T}'$ we
obtain a subfamily $\mathcal{T}''$, being the union of
$\mathcal{T}'$, with ${\rm RS}(\mathcal{T}')=2$ and ${\rm
ds}(\mathcal{T}')=n$.

Now we can continue the process for greater countable ordinals
$\alpha$ obtaining a $d_\infty$-definable subfamily
$\mathcal{T}^\ast\subset\mathcal{T}$ with ${\rm
RS}(\mathcal{T}^\ast)=\alpha$ and ${\rm ds}(\mathcal{T}^\ast)=n$
for given $n\in\omega\setminus\{0\}$.~$\Box$

\medskip
Theorem 4.6 and Remark 4.7 imply the following:

\medskip
{\bf Theorem 4.8.} {\em Let $\mathcal{T}$ be a family of a
countable language $\Sigma$ and with ${\rm
RS}(\mathcal{T})=\infty$, $\alpha$ be a countable ordinal,
$n\in\omega\setminus\{0\}$. Then there is a $d_\infty$-definable
subfamily $\mathcal{T}^\ast\subset\mathcal{T}$ such that ${\rm
RS}(\mathcal{T}^\ast)=\alpha$ and ${\rm ds}(\mathcal{T}^\ast)=n$.}

\medskip
{\bf Example 4.9.} Let $\mathcal{T}_\Sigma$ be a family of all
theories of a countable language $\Sigma$ with ${\rm
RS}(\mathcal{T}^\ast)=\infty$ \cite{MarkhSud}, say of unary
predicates $Q_n$, $n\in\omega$. Taking a countable $d$-definable
subfamily $\mathcal{T}\subset\mathcal{T}_\Sigma$ with either empty
or complete predicates $Q_n$ such that complete predicates in
$\mathcal{T}$ are linearly ordered and indexes for complete
predicates form an infinite set $I\subset\omega$ with infinite
$\omega\setminus I$ we can assume that $\mathcal{T}$ is
$e$-minimal has unique accumulation point witnessing ${\rm
RS}(\mathcal{T})=1$ and ${\rm ds}(\mathcal{T})=1$. Taking indexes
in $\omega\setminus I$ we can define countably many disjoint
$e$-minimal $d$-definable subfamilies $\mathcal{T}_k$ with unique
accumulation point for the set of all accumulation points of
$\mathcal{T}_k$ witnessing ${\rm RS}=2$ and ${\rm ds}=1$. Now
applying Theorem 4.8 we can continue the process obtaining ${\rm
RS}=\alpha$ and ${\rm ds}=n$ for arbitrary countable ordinal
$\alpha$ and $n\in\omega\setminus\{0\}$.~$\Box$

\medskip
{\bf Definition.} An $\alpha$-ranking set $\Phi$ for
$\mathcal{T}$, and $\mathcal{T}_\Phi$ are called {\em
$\mathcal{T}$-irreducible} if for any $\mathcal{T}$-inconsistent
$\Psi,{\rm X}\supseteq\Phi$, i.e., with
$\mathcal{T}_\Psi\cap\mathcal{T}_{\rm
X}\supseteq\mathcal{T}_\Phi$, ${\rm RS}(\mathcal{T}_\Psi)<\alpha$
or ${\rm RS}(\mathcal{T}_{\rm X})<\alpha$. An $\alpha$-ranking
sentence $\varphi$ for $\mathcal{T}$, and $\mathcal{T}_\varphi$
are called {\em $\mathcal{T}$-irreducible} if the singleton
$\{\varphi\}$ is $\mathcal{T}$-irreducible.

If $\mathcal{T}$ is fixed, $\mathcal{T}$-irreducible sets are
called simply {\em irreducible}.

\medskip
By the definition each $\mathcal{T}$-inconsistent set $\Phi$, with
$\mathcal{T}_\Phi=\emptyset$, is irreducible, as well as
singletons $\mathcal{T}_\Phi$.

Moreover, nonempty $E$-closed families $\mathcal{T}_\Phi$ of rank
$\alpha$ are irreducible if and only if ${\rm
ds}(\mathcal{T}_\Phi)=1$.

Indeed, if $\mathcal{T}_\Phi$ is irreducible it can not be divided
by a sentence into two parts of rank $\alpha$ implying ${\rm
ds}(\mathcal{T}_\Phi)=1$. Conversely, having
$\mathcal{T}$-inconsistent $\Psi,{\rm X}\supseteq\Phi$, with
$\mathcal{T}_\Psi\cap\mathcal{T}_{\rm
X}\supseteq\mathcal{T}_\Phi$, ${\rm RS}(\mathcal{T}_\Psi)=\alpha$
and ${\rm RS}(\mathcal{T}_{\rm X})=\alpha$, we obtain, by
compactness, some $\mathcal{T}$-inconsistent $\psi\in\Psi$ and
$\chi\in{\rm X}$ such that ${\rm
RS}((\mathcal{T}_\Phi)_\psi)=\alpha$ and ${\rm
RS}((\mathcal{T}_\Phi)_\chi)=\alpha$ contradicting ${\rm
ds}(\mathcal{T}_\Phi)=1$.

Since each family $\mathcal{T}$ with ${\rm
RS}(\mathcal{T})=\alpha\geq 0$ has a finite degree ${\rm
ds}(\mathcal{T})=n$, there are pairwise inconsistent sentences
$\varphi_1,\ldots,\varphi_n$ such that
$\mathcal{T}=\mathcal{T}_{\varphi_1}\,\dot{\cup}\,\ldots$
$\,\dot{\cup}\,\mathcal{T}_{\varphi_n}$, ${\rm
RS}(\mathcal{T}_{\varphi_i})=\alpha$ and ${\rm
ds}(\mathcal{T}_{\varphi_i})=1$, $i=1,\ldots,n$.

Thus, all $e$-totally transcendental $E$-closed families and, in
particular, $d$-definable $\alpha$-ranking $E$-closed families are
reduced to irreducible ones:

\medskip
{\bf Proposition 4.10.} {\em Any $e$-totally transcendental
$E$-closed family $\mathcal{T}$ is represented as a finite
disjoint union of $s$-definable irreducible subfamilies of rank
$\alpha={\rm RS}(\mathcal{T})$.}


\end{document}